\documentclass[journal]{elsarticle}

\usepackage{lineno,hyperref}
\usepackage{epsfig,extarrows}
\usepackage{amsmath,amsthm,subfigure}
\usepackage{amssymb}
\usepackage{graphicx,epstopdf,latexsym,amsfonts,fancyhdr}
\usepackage{booktabs,arydshln,multirow}
\usepackage{setspace}
\usepackage{xcolor}
\usepackage{stmaryrd}
\usepackage[utf8]{inputenc}
\usepackage{lineno,hyperref}
\usepackage[top=1in, bottom=1in, left=1in, right=1in]{geometry}

\newtheorem{theorem}{Theorem}[section]

\newtheorem{lemma}[theorem]{Lemma}

\newtheorem{definition}[theorem]{Definition}
\newtheorem{remark}[theorem]{Remark}

\modulolinenumbers[1]

\journal{Journal of \LaTeX\ Templates}

\begin{document}

\begin{frontmatter}

\title{A Pick function approach for designing energy-decay preserving schemes of the Maxwell equations in Havriliak-Negami dispersive media}

\author[mymainaddress,mysecondaddress]{Baoli Yin\corref{mycorrespondingauthor}}
\cortext[mycorrespondingauthor]{Corresponding author}
\ead{baolimath@126.com}

\author[mymainaddress,mysecondaddress]{Guoyu Zhang}
\ead{guoyu\_zhang@imu.edu.cn}

\author[mymainaddress,mysecondaddress]{Yang Liu}
\ead{mathliuyang@imu.edu.cn}

\author[mymainaddress,mysecondaddress]{Hong Li}
\ead{smslh@imu.edu.cn}

\address[mymainaddress]{School of Mathematical Sciences, Inner Mongolia University, Hohhot 010021, China;}
\address[mysecondaddress]{Inner Mongolia Key Laboratory of Mathematical Modeling and Scientific Computing, Hohhot 010021, China;}

\begin{abstract}
This work proposes a novel approach for designing high-order energy-decaying schemes for Maxwell’s equations in Havriliak–Negami dispersive media. It is shown that conventional convolution quadrature (CQ) methods, which rely directly on the generating function of linear multistep methods, cannot generate completely monotonic sequences beyond first-order accuracy. We rigorously prove that for any linear multistep method of second- or higher-order, the associated generating function \(\delta(\zeta)\) cannot satisfy both that \(-\delta(\zeta)\) is a Pick function and that it is analytic on \((-\infty,1)\) — a key requirement for constructing completely monotonic sequences.
To overcome this fundamental limitation, we introduce a reconstruction of the generating function’s structure. By strategically incorporating the theory of Pick functions, we successfully construct a second-order completely monotonic sequence. This theoretical advance leads to a discrete scheme that inherits the continuous model’s energy decay property, guaranteeing unconditional stability. Numerical experiments confirm the convergence rates and energy dissipation behavior of the proposed method.
\end{abstract}


\end{frontmatter}


\section{Intruction}
Building upon its ability to describe asymmetrical and broad relaxation spectra, the Havriliak-Negami (H-N) model \cite{havriliak1966complex} has become a critical tool in the dielectric characterization of diverse materials systems. 
Its empirical form is extensively applied across multiple fields, including: Polymer Science and Soft Matter Physics, where it characterizes dielectric and mechanical relaxations \cite{havriliak1967complex}, dynamic glass transitions \cite{mckenna201750th}, and dynamics in microporous polymer membranes \cite{rodriguez2024penetrant}; Condensed Matter Physics, for describing relaxation and diffusion in disordered systems such as amorphous semiconductors \cite{ngai2023universal,elliott1987ac}; Geophysics, in interpreting dielectric spectra of sedimentary rocks \cite{lesmes2001dielectric} and subsurface characterization \cite{hubbard2011hydrogeophysics}; Liquid Crystal Research, for analyzing dielectric behavior in ferroelectric nematic phases \cite{sebastian2020ferroelectric,mandle2021molecular}; and Biophysics, where its link to fractional calculus supports the modeling of relaxation in biological tissues \cite{ionescu2017role}.
\par
The H-N relation's generality is demonstrated by its reduction to several classical models under specific parameter constraints: it simplifies to the Debye model \cite{debye1929polar} when \(\alpha = \beta = 1\), to the Davidson-Cole (D-C) model \cite{davidson1951dielectric} for \(\alpha = 1\), and to the Cole-Cole (C-C) model \cite{rh1941dispersion} for \(\beta = 1\). This flexibility allows it to unify the description of anomalously dispersive dielectric behavior.
This widespread utility stems from the model's comprehensive mathematical framework, which generalizes the dielectric response of dispersive media. The complex relative permittivity \(\epsilon_r(\omega)\) is defined by the expression:
\[
\epsilon_r(\omega) = \epsilon_\infty + \frac{\epsilon_s - \epsilon_\infty}{(1 + (i\omega\tau_0)^\alpha)^\beta},
\]
where \(0 < \alpha, \beta \leq 1\) are shape parameters. The other key parameters are the high-frequency limit \(\epsilon_\infty\), the static permittivity \(\epsilon_s\) (with \(\epsilon_s > \epsilon_\infty \geq 1\)), the characteristic relaxation time \(\tau_0\), the angular frequency \(\omega\), and the imaginary unit \(i\).
After rescaling the variables \cite{yang2021analysis},  the time domain Maxwell's equations in Havriliak-Negami(H-N) dispersive medium can be formulated by
\begin{align}
    \epsilon_{\infty} \partial_t \boldsymbol E + \partial_t \boldsymbol P &= \nabla \times \boldsymbol H, \quad (\boldsymbol x,t) \in \Omega \times (0,T], \tag{1.1}\label{a1} \\
    \partial_t \boldsymbol H &= -\nabla \times \boldsymbol E, \quad (\boldsymbol x,t) \in \Omega \times (0,T], \tag{1.2}\label{a2} \\
    \boldsymbol P(x,t) &= \Delta \epsilon \int_0^t \omega_{\alpha,\beta}(t-s)\boldsymbol E(x,s)\,ds, \quad (\boldsymbol x,t) \in \Omega \times (0,T], \tag{1.3}\label{a3}
\end{align}
with the following initial conditions:
\begin{align}
    \boldsymbol E(x,0) &=\boldsymbol E_0(x), \quad \boldsymbol H(x,0) = \boldsymbol H_0(x), \quad \boldsymbol P(x,0) = 0, \text{ for } \boldsymbol x \in \Omega, \tag{1.4}\label{a4}
\end{align}
and perfect conducting boundary condition:
\begin{align}
   \boldsymbol n \times \boldsymbol E &= 0 \quad \text{ on } \quad \partial \Omega \times (0,T). \tag{1.5}\label{a5}
\end{align}
where $\Delta\epsilon:=\epsilon_s-\epsilon_\infty$ and the kernel $\omega_{\alpha,\beta}(t)$ is defined by the inverse Laplace transform: $\mathcal{L}^{-1}\big[(1+s^\alpha)^{-\beta}\big](t)$.
Generally, $\omega_{\alpha,\beta}(t)$ cannot be expressed in a closed form.
In fact, $\omega_{\alpha,\beta}(t)$ can be expressed by resorting to the Prabhakar function \cite{prabhakar1971singular}, which is a Mittag-Leffler function with three parameters.
\par
The numerical investigation of the Cole-Cole (C-C) model has been extensive due to its relative simplicity. For instance, Li et al. \cite{li2012time} established that the continuous energy does not exceed the initial energy, while Yin et al. \cite{baoli2023discrete} subsequently proposed a criterion for assessing the energy stability of time-stepping schemes. Additional related contributions are documented in \cite{huang2019accurate, xiao2025unconditionally, xie2022efficient}. Furthermore, in \cite{zhang2025analysis}, a second-order energy-decaying $\theta$-scheme (with $\theta=\frac{1}{2}$) was developed, ensuring the discrete energy satisfies \(\mathbb{E}^n \leq \mathbb{E}^{n-1}\).
In contrast, numerical studies for the more general Havriliak-Negami (H-N) model have predominantly focused on frequency-domain simulations using Debye-model approximations \cite{chakarothai2018novel, kelley2007debye, rekanos2012fdtd}. These approaches, however, generally lack rigorous stability or convergence analysis. Yang et al. \cite{yang2021analysis} proposed an energy-stable scheme for the H-N model, demonstrating that the discrete energy at any time remains bounded by the initial energy. They further showed that a modified discrete energy incorporating a history term satisfies a decay property. Nevertheless, their method is limited to first-order accuracy, and a continuous energy-decay law for the H-N model remained unestablished.
\par
In this work, we first derive a continuous energy-decay property for the H-N model. We then present, for the first time, a second-order accurate numerical scheme that preserves a discrete energy-decay property. This is achieved by constructing a complete monotonicity-preserving approximation for the convolution term \(\omega_{\alpha,\beta} * \boldsymbol{E}\), which combines the theory of Pick functions \cite{donoghue2012monotone} with convolution quadrature (CQ) techniques \cite{lubich1988convolution}.
Specifically, we depart from the structure of the generating function proposed in CQ,  
\[
w(\zeta) = \left(1 + \left( \frac{\delta(\zeta)}{\tau} \right)^\alpha \right)^{-\beta},
\]  
and instead introduce the following form:  
\[
w(\zeta)=\bigg(1+\bigg(\frac{1-\zeta}{\tau}\bigg)^\alpha G(\zeta)^{1-\alpha}\bigg)^{-\beta},
\]  
under the assumption that \(-G(\zeta)\) is a Pick function.  
This modified generating function yields a completely monotonic sequence \((w_0, w_1, \dots )\), and the discrete convolution \(\displaystyle\sum_{k=0}^n w_{n-k}\boldsymbol E^k\) achieves second-order accuracy in approximating \(\omega_{\alpha,\beta} * \boldsymbol{E}\).
For further background on complete monotonicity and its advantageous properties, we refer to \cite{dong2025mittag, li2021complete} and references therein.
\par
In summary, our main contributions are:
\begin{itemize}
\item Deriving a continuous energy-decay law for the H-N model based on the complete monotonicity of its kernel function.
\item Rigorously proving that no second-order linear multistep method capable of generating completely monotonic sequences can directly yield a second-order scheme within the CQ framework.
\item Introducing a novel second-order, complete-monotonicity-preserving numerical formula by synergistically employing Pick functions and convolution quadrature theory.
\end{itemize}
\par
This paper is organized as follows.
Section 2 introduces the necessary mathematical preliminaries, including notations and key properties related to Prabhakar functions, Pick functions, and completely monotonic sequences.
In Section 3, we rigorously derive a continuous energy-decay law for the Havriliak-Negami model by leveraging the complete monotonicity of its underlying kernel function.
Section 4 first reviews relevant concepts from convolution quadrature theory, followed by a proof that no second-order linear multistep methods capable of generating completely monotonic sequences that can directly produce a second-order scheme for this problem. 
We then present a novel second-order, complete-monotonicity-preserving approximation by utilizing the theory of Pick functions, enabling the construction of discrete energy-decay-preserving schemes.
Section 5 presents several numerical experiments to validate the theoretical findings.
Concluding remarks are provided in Section 6.
\par
We first fix the notation and functional setting used throughout this work.
For \( r \geq 0 \), let \( H^r(\Omega) \) denote the standard Sobolev space equipped with norm \( \|\cdot\|_r \). In particular, we write \( \|\cdot\| \) for \( \|\cdot\|_0 \), and note that \( H^0(\Omega) \) coincides with \( L^2(\Omega) \).
Define the curl-conforming Sobolev space
\[
H^r(\operatorname{curl};\Omega) = \left\{ \boldsymbol{v} \in \bigl(H^r(\Omega)\bigr)^2 : \nabla \times \boldsymbol{v} \in H^r(\Omega) \right\},
\]
endowed with the norm
\[
\|\boldsymbol{v}\|_{r,\operatorname{curl}} = \left( \|\boldsymbol{v}\|_r^2 + \|\nabla \times \boldsymbol{v}\|_r^2 \right)^{1/2}.
\]
Its subspace with vanishing tangential trace is given by
\[
H_0(\operatorname{curl};\Omega) = \left\{ \boldsymbol{v} \in H(\operatorname{curl};\Omega) : \boldsymbol{n} \times \boldsymbol{v} = \mathbf{0} \text{ on } \partial\Omega \right\},
\]
where we denote \( H(\operatorname{curl};\Omega) = H^0(\operatorname{curl};\Omega) \).
\section{Preliminaries}
\subsection{Prabhakar Function}
The Mittag-Leffler function with three parameters, i.e., the Prabhakar function \cite{prabhakar1971singular}, is defined by
\begin{equation}\label{a6}\begin{split}
\mathcal{E}^\gamma_{\rho,\mu}(z)=\frac{1}{\Gamma(\gamma)}\sum_{k=0}^{\infty}\frac{\Gamma(k+\gamma)}{\Gamma(\rho k+\mu)}\frac{z^k}{k!},\quad \Re (\rho)>0, \Re (\mu)>0, \gamma>0.
\end{split}
\end{equation}
Let $e_{\rho,\mu}^\gamma(t;\lambda)
=t^{\mu-1}\mathcal{E}^\gamma_{\rho,\mu}(\lambda t^\rho)$, there holds \cite[p.47]{kilbas2006theory},
\begin{equation}\label{a7}\begin{split}
\mathcal{L}\big[e_{\rho,\mu}^\gamma(t;\lambda)\big](s)
=\frac{s^{\rho\gamma-\mu}}{(s^\rho-\lambda)^\gamma},
\end{split}
\end{equation}
which means $\omega_{\alpha,\beta}(t)=\mathcal{L}^{-1}\big[(1+s^\alpha)^{-\beta}\big](t)
=e_{\alpha,\alpha\beta}^\beta(t;-1)$ by using the substitution $\rho=\alpha, \lambda=-1, \gamma=\beta, \mu=\alpha\beta$.
The kernel function $e_{\alpha,\alpha\beta}^\beta(t;-\lambda)$ is related to the Prabhakar fractional integral of order $\alpha,\beta>0$ with $\lambda>0$ defined by
\begin{equation}\label{b1.1}
(\mathcal{J}_{0,t}^\alpha+\lambda)^\beta g(t)=\int_{0}^{t}e_{\alpha,\alpha\beta}^\beta(t-s;-\lambda)g(s)\mathrm{d}s.
\end{equation}
It is known that  for $g(t)=t^k$ with $k>-1$,
\begin{equation}\label{b1.2}
(\mathcal{J}_{0,t}^\alpha+\lambda)^\beta g(t)=\Gamma(k+1) e_{\alpha,\alpha\beta+k+1}^\beta(t;-\lambda).
\end{equation}
\begin{definition}[Completely monotonic functions]
A function $\omega:(0,\infty) \to \mathbb{R}$ beloning to $C^\infty$ is called completely monotonic if
\begin{equation}\begin{split}\label{b000}
(-1)^n \omega^{(n)}(t) \geq 0, \quad \forall t >0, n=0,1,2,\cdots.
\end{split}\end{equation}
\end{definition}
By the Bernsterin's theorem \cite{bernstein1929fonctions, widder2015laplace, feller1991introduction}, 
a function $f$ is completely monotonic if and only if there exists a positive measure $\mu(\lambda)$ such that $\omega(t)=\int_0^\infty e^{-t\lambda}\mathrm{d}\mu(\lambda)$.

The function $e_{\rho,\mu}^\gamma(t;\lambda)$ is completely monotonic \cite{mainardi2015complete,de2011models} for $0<\rho\leq 1, 0<\rho\gamma \leq \mu \leq 1$, indicating that $\omega_{\alpha,\beta}(t)$ of the H-N model is also completely monotonic.

\subsection{Pick Functions}
Pick functions are key tools in our study to construct high-order energy-decay numerical schemes for the H-N model.
Denote by $\mathbb{C}$ the complex plain and by $\mathbb{C}^+$  the open upper half-plane.
\begin{definition}\cite{donoghue2012monotone}
A function $\phi(\zeta)=U(\zeta)+{\rm i}V(\zeta)$ for $\zeta\in \mathbb{C}$ is called a Pick/Nevanlinna function if $\phi$ is analytic in  $\mathbb{C}^+$, and has non-negative imaginary part, i.e., $V(\zeta)\geq 0$.
\end{definition}
Note that if a Pick function attains a real value at any point in $\mathbb{C}^+$ , it must be identically constant, according to the maximum modulus principle. To exclude this trivial case, we therefore assume in the following analysis that \textbf{the Pick function under consideration is nonconstant}, i.e., $V(\zeta) \not\equiv 0$.
Given this assumption, it is straightforward to verify that the composition of two Pick functions again yields a Pick function.
\begin{lemma}\cite{donoghue2012monotone}\label{lem.3}
      Any rational Pick function $\phi(\zeta)$ which is real on an interval of the real axis must take exactly the following form:
      \begin{equation}
        \phi(\zeta)=\gamma_1 \zeta+\gamma_0+\sum_{i=0}^{N}\frac{m_i}{\lambda_i-\zeta},
        \quad
        \text{where } \gamma_1\geq0, m_i>0 ~\text{ and $\gamma_0$ is real}.
      \end{equation}
  \end{lemma}

\subsection{Completely Monotonic Sequences}
We say a sequence $(v_0,v_1,\cdots)$ is completely monotonic \cite{liu2016generating}, or has complete monotonicity property, if
\begin{equation}\label{b00}
(I-S)^kv_j=\sum_{n=0}^{k}(-1)^n
\begin{pmatrix}
  k \\
  n
\end{pmatrix}
v_{n+j}\geq 0,\quad k=0,1,2,\cdots,
\end{equation}
where $S$ is the backshift operator: $Sv_j=v_{j+1}$ for $j\geq 0$.
 \begin{lemma}\label{lem.2}\cite{liu2016generating}
    A sequence $(v_0,v_1,\cdots)$ is completely monotonic if and only if its generating function $v(\zeta)=\displaystyle\sum_{j=0}^{\infty}v_j\zeta^j$ is a Pick function that is analytic and nonnegative on $(-\infty,1)$.
  \end{lemma}
\section{Energy-decay law of H-N model}
In \cite{yang2021analysis}, the authors presented the following energy stable property of the H-N model (\ref{a1})-(\ref{a5})
\begin{equation}\begin{split}
\epsilon_\infty \|\boldsymbol{E}(\cdot,t)\|^2+\|\boldsymbol{H}(\cdot,t)\|^2 
\leq
\epsilon_\infty \|\boldsymbol{E}_0\|^2+\|\boldsymbol{H}_0\|^2,\quad \forall t \in (0,T].
\end{split}\end{equation}
This section is devoted to establishing rigorously the monotonic decay of the continuous energy associated with the H-N model.
\begin{lemma}\label{lem.1}
Assume $u(t)$ is continuous and bounded on $[0,\infty)$ and $\omega(t)$ is completely monotonic.
For any $t>0$, there holds
\[
u(t) \frac{d}{dt} \int_0^t \omega(t-s) u(s)  \mathrm{d}s \geq \frac{1}{2} \frac{d}{dt} \int_0^t \omega(t-s) u^2(s) \mathrm{d}s.
\] 
\begin{proof}
Since $\omega(t)$ is a completely monotonic function, by Bernstein’s theorem, there exists a positive measure $\mu$ on $[0,\infty)$ such that
$$
\omega(t) = \int_0^\infty e^{-\lambda t}  d\mu(\lambda), \quad t > 0.
$$
Define the convolution integrals
$$
A(t) = \int_0^t \omega(t-s) u(s)  ds, \quad B(t) = \int_0^t \omega(t-s) u^2(s)  ds.
$$
We aim to prove the inequality
$$
u(t) \frac{d}{dt} A(t) \geq \frac{1}{2} \frac{d}{dt} B(t).
$$
Using the Laplace representation of $\omega$, we write
$$
A(t) = \int_0^\infty e^{-\lambda t} U_\lambda(t)  d\mu(\lambda), \quad B(t) = \int_0^\infty e^{-\lambda t} V_\lambda(t)  d\mu(\lambda),
$$
where
$$
U_\lambda(t) = \int_0^t e^{\lambda s} u(s)  ds, \quad V_\lambda(t) = \int_0^t e^{\lambda s} u^2(s)  ds.
$$
Differentiating under the integral sign yields
$$
\frac{d}{dt} A(t) = \int_0^\infty \left[ -\lambda e^{-\lambda t} U_\lambda(t) + u(t) \right] d\mu(\lambda),
\quad
\frac{d}{dt} B(t) = \int_0^\infty \left[ -\lambda e^{-\lambda t} V_\lambda(t) + u^2(t) \right] d\mu(\lambda).
$$
Thus,
$$
D(t): = u(t) \frac{d}{dt} A(t) - \frac{1}{2} \frac{d}{dt} B(t)
= \int_0^\infty \left[ -\lambda e^{-\lambda t} U_\lambda(t) u(t) + \frac{1}{2} u^2(t) + \frac{1}{2} \lambda e^{-\lambda t} V_\lambda(t) \right] d\mu(\lambda).
$$
Now consider the identity
$$
\lambda e^{-\lambda t} V_\lambda(t) = 2 \lambda e^{-\lambda t} u(t) U_\lambda(t) - u^2(t) (1 - e^{-\lambda t}) + \lambda e^{-\lambda t} \int_0^t e^{\lambda s} [u(t) - u(s)]^2  ds,
$$
which follows from expanding $\int_0^t e^{\lambda s} [u(t) - u(s)]^2  ds$. Substituting into $D(t)$ and simplifying gives
$$
D(t) = \int_0^\infty \left[ \frac{1}{2} u^2(t) e^{-\lambda t} + \frac{1}{2} \lambda e^{-\lambda t} \int_0^t e^{\lambda s} [u(t) - u(s)]^2  ds \right] d\mu(\lambda).
$$
Since $\int_0^\infty e^{-\lambda t}  d\mu(\lambda) = \omega(t)$, we obtain
$$
D(t) = \frac{1}{2} u^2(t) \omega(t) + \frac{1}{2} \int_0^\infty \lambda e^{-\lambda t} \left( \int_0^t e^{\lambda s} [u(t) - u(s)]^2  ds \right) d\mu(\lambda).
$$
As $\mu$ is a positive measure, $\lambda \geq 0$, and the integrands are nonnegative, it follows that $D(t) \geq 0$
which completes the proof.
\begin{remark}
We note that Lemma \ref{lem.1} extends the findings in \cite{alsaedi2015maximum} where the Riemann-Liouville fractional operator is involved with the kernel function $\omega(t)=\frac{1}{\Gamma(1-\alpha)}\frac{1}{t^\alpha}$.
It is straightforward to verify that $t^{-\alpha}$ for $0<\alpha<1$ is, in fact, a completely monotonic function.
\end{remark}
\end{proof} 
\end{lemma}
\begin{theorem}
Assume $\boldsymbol{E}_0, \boldsymbol{H}_0 \in \boldsymbol{L}^2(\Omega)$.
The H-N model (\ref{a1})-(\ref{a5}) satisfies the energy decay property
\begin{equation}\begin{split}\label{b0}
\mathbb{E}(t_1) \leq \mathbb{E}(t_2), \quad \forall t_1,t_2 \in [0,T], t_1>t_2,
\end{split}\end{equation}
where the energy $\mathbb{E}(t)$ is defined by
\begin{equation}\begin{split}
\mathbb{E}(t)=\epsilon_\infty \|\boldsymbol{E}(\cdot,t)\|^2+\|\boldsymbol{H}(\cdot,t)\|^2 
+\Delta\epsilon\int_{0}^{t}\omega_{\alpha,\beta}(t-s) \|\boldsymbol E(\cdot,s)\|^2\mathrm{d}s.
\end{split}\end{equation}
\begin{proof}
We first give the weak form of the H-N model, i.e., find \( \boldsymbol E \in H_0(\mathrm{curl}, \Omega) \)  and  \( \boldsymbol H, \boldsymbol P \in \boldsymbol L^2(\Omega) \)  so that
\begin{align}
    &\epsilon_\infty (\partial_t \boldsymbol E, \boldsymbol\phi) + (\partial_t \boldsymbol P, \boldsymbol\phi) - (\boldsymbol H, \nabla \times \boldsymbol\phi) = 0, 
    & \forall \boldsymbol\phi \in H_0(\mathrm{curl}, \Omega), \tag{2.1}\label{b1} \\
    &(\partial_t \boldsymbol H, \boldsymbol\psi) + (\nabla \times \boldsymbol E, \boldsymbol\psi) = 0, 
    & \forall \boldsymbol\psi \in \boldsymbol L^2(\Omega), \tag{2.2}\label{b2}  \\
    &(\boldsymbol P, \boldsymbol\varphi) - \Delta \epsilon \int_0^t \omega_{\alpha, \beta} (t - s) \, (\boldsymbol E(\cdot, s), \boldsymbol\varphi) \, ds=0, 
    & \forall \boldsymbol\varphi \in \boldsymbol L^2(\Omega). \tag{2.3}\label{b3} 
\end{align}
We begin by setting $\boldsymbol\phi = \boldsymbol E$ in equation (\ref{b1}) and $\boldsymbol\psi = \boldsymbol H$ in equation (\ref{b2}), then adding the two resulting equations to yield
\[
\varepsilon_\infty \left( \partial_t \boldsymbol E, \boldsymbol E \right) + \left(\partial_t \boldsymbol H, \boldsymbol H \right) + (\partial_t\boldsymbol P, \boldsymbol E) = 0.
\]
Replacing $\boldsymbol P$ by equation (\ref{a3}) and using the fact that $\omega_{\alpha,\beta}(t)$ is a completely monotonic function, one immediately obtain (\ref{b0})  in accordance to Lemma \ref{lem.1}.
\end{proof}
\end{theorem}
\section{Discrete energy decay schemes}
\subsection{Recall the CQ Theory}
We begin by recalling the fundamentals of convolution quadrature theory \cite{lubich1988convolution}. To approximate the convolution \( \int_0^t \omega(t-s)u(s)\, ds \) at a discrete time \( t_n \), the CQ method employs a discrete sum \( \displaystyle\sum_{k=0}^n w_{n-k}(\tau)u^k \), where \( \tau \) denotes the time step size.
The sequence $(w_0, w_1, \cdots )$ is generated by the generating function $w(\zeta)$:
\begin{equation}\begin{split}
\sum_{k=0}^\infty w_k \zeta^k= w(\zeta):=W\bigg(\frac{\delta(\zeta)}{\tau}\bigg),
\quad
W(s)=\mathcal{L}[\omega(t)]=(1+s^\alpha)^{-\beta},
\end{split}\end{equation}
where $\delta(\zeta)$ is the quotient of the generating polynomials of a linear multistep method when solving the differential equation $y'=\lambda y+u$.
Specifically, consider the following linear multistep method
\begin{equation}\begin{split}
\sum_{j=0}^\ell a_j y^{n+j-\ell}=\tau\sum_{j=0}^\ell b_j (\lambda y^{n+j-\ell} +u(t_{n+j-\ell})),\quad n \geq 0,
\end{split}\end{equation}
where $y_{-\ell}=\dots=y_{-1}=0$ and $u \in C[0,\infty)$ is extened by zero for $t<0$.
Then $\delta(\zeta)$ is defined by
\begin{equation}\begin{split}\label{b11}
\delta(\zeta)=\frac{a_0\zeta^\ell+a_1\zeta^{\ell-1}+\cdots+a_\ell}{b_0\zeta^\ell+b_1\zeta^{\ell-1}+\cdots+b_\ell}.
\end{split}\end{equation}
A well know fact is that a linear multistep method is consistent of order $p$ if and only if
\begin{equation}\begin{split}\label{b6}
\frac{1}{\tau}\delta(e^{-\tau})=1+O(\tau^p).
\end{split}\end{equation}
Examples are the backward difference formulas of order $p$ (BDF-$p$):
$\displaystyle\delta(\zeta)=\sum_{k=1}^p \frac{1}{k}(1-\zeta)^k$ for $p=1,2,\cdots,6$.
\par
Accordingly, for the H-N model at time \( t_n \), the convolution integral \( \int_0^{t_n} \omega_{\alpha,\beta}(t_n-s)\boldsymbol{E}(x,s)\,ds \) is discretized via the CQ method. This yields the approximation:
\begin{equation}
\begin{split}\label{b4}
\sum_{k=0}^n w_{n-k}(\tau)\boldsymbol{E}^k, \quad \text{where} \quad
w(\zeta) = \left(1 + \left( \frac{\delta(\zeta)}{\tau} \right)^\alpha \right)^{-\beta}.
\end{split}
\end{equation}
This discretization possesses an asymptotic accuracy of order \( p \), as ensured by the underlying CQ theory. Specifically, the approximation error satisfies:
\begin{equation}
\begin{split}
\left\| \sum_{k=0}^n w_{n-k}(\tau)\boldsymbol{E}^k - \int_0^{t_n} \omega_{\alpha,\beta}(t_n-s)\boldsymbol{E}(x,s),ds \right\| = O(\tau^p).
\end{split}
\end{equation}
\subsection{Construction of a Completely Monotonic Sequence}
Given that the kernel function $\omega_{\alpha,\beta}(t)$ is completely monotonic, it is imperative that the discrete sequence $(w_0, w_1, \cdots)$ also exhibits complete monotonicity. According to Lemma \ref{lem.2}, a necessary condition for this is that the generating function $w(\zeta)$ must be a Pick function.
\begin{lemma}
If \( -\delta(\zeta) \) is a Pick function, then \( w(\zeta) \), as defined in (\ref{b4}), is also a Pick function.
\begin{proof}
Since \(-\delta(\zeta)\) is a Pick function, \(\delta(\zeta)\) is analytic on \(\mathbb{C}^+\) and satisfies \(\text{Im}(\delta(\zeta)) < 0\) for all \(\zeta \in \mathbb{C}^+\). 
Consider the function \(W(z) = (1 + z^\alpha)^{-\beta}\), where complex powers are defined using the principal branch. 
For \(z\) with \(\text{Im}(z) < 0\), we have \(\arg(z) \in (-\pi, 0)\), so \(\arg(z^\alpha) = \alpha\arg(z) \in (-\alpha\pi, 0) \subset (-\pi, 0)\). 
Thus, \(z^\alpha\) remains in the lower half-plane with \(\text{Im}(z^\alpha) < 0\).

Now, \(1 + z^\alpha\) also lies strictly in the lower half-plane since \(\text{Im}(1 + z^\alpha) = \text{Im}(z^\alpha) < 0\), and it avoids the negative real axis due to the strict inequality. Therefore, \(\arg(1 + z^\alpha) \in (-\pi, 0)\), and we obtain:
\[
\arg\left((1 + z^\alpha)^{-\beta}\right) = -\beta\arg(1 + z^\alpha) \in (0, \beta\pi) \subset (0, \pi),
\]
which implies \(\text{Im}(W(z)) > 0\). Hence, \(W(z)\) maps the lower half-plane strictly to the upper half-plane.

Composing with \(z = \delta(\zeta)/\tau\), which maps \(\mathbb{C}^+\) to the lower half-plane, yields \(w(\zeta) = W(\delta(\zeta)/\tau)\). This composition is analytic on \(\mathbb{C}^+\) and satisfies \(\text{Im}(w(\zeta)) > 0\) for \(\zeta \in \mathbb{C}^+\), confirming that \(w(\zeta)\) is a Pick function.
\end{proof}
\end{lemma}
While the function $-\delta(\zeta)=\zeta-1$ for the Euler method (BDF-1) is a Pick function and is analytic for $\zeta \in (-\infty,1)$, this property does not extend to the $\delta(\zeta)$ functions of higher-order linear multistep methods.

\begin{lemma}
No second-order linear multistep method exists for which $-\delta(\zeta)$ is a Pick function and is analytic for $\zeta \in (-\infty,1)$.
\begin{proof}
Assume $-\delta(\zeta)$ is a Pick function.
Since $\delta(\zeta)$ is a rational function according to (\ref{b11}), by Lemma \ref{lem.3}, $\delta(\zeta)$ must be of the form
\begin{equation}\begin{split}
\delta(\zeta)=
-\gamma_1 \zeta-\gamma_0-\sum_{i=0}^{N}\frac{m_i}{\lambda_i-\zeta},
        \quad
        \text{where } \gamma_1\geq0, m_i>0 ~\text{ and $\gamma_0$ is real}.
\end{split}\end{equation}
If $-\delta(\zeta)$ is further analytic at $(-\infty,1)$, then $\lambda_i \geq 1$.
\par
We first show that $\lambda_i \neq 1, \forall i$. 
This is for the reason that $\frac{m_i}{1-\zeta}|_{\zeta=e^{-\tau}}$ tends to $+\infty$ as $\tau \to 0$, which contradicts the fact that $\delta(e^{-\tau}) \to 0$ by the condition (\ref{b6}).

Assume $\lambda_i>1$. Let $\zeta=e^{-\tau}=1-\tau+\frac{\tau^2}{2}+R_1(\tau)$ where $R_1(\tau)=O(\tau^3)$.
Then,
\begin{equation}\begin{split}
\frac{m_i}{\lambda_i-\zeta}\bigg|_{\zeta=e^{-\tau}}
&=\frac{m_i}{\lambda_i-1}
\sum_{k=0}^2 \bigg(
\frac{e^{-\tau}-1}{\lambda_i-1}
\bigg)^k
+O(\tau^3)
\\&=
\frac{m_i}{\lambda_i-1}
\bigg[
1+\frac{\tau}{1-\lambda_i}
+\frac{1+\lambda_i}{2(1-\lambda_i)^2}\tau^2
\bigg]
+O(\tau^3),
\end{split}\end{equation}
yielding that
\begin{equation}\begin{split}
\delta(e^{-\tau})=
\delta_0+\delta_1\tau+\delta_2\tau^2+O(\tau^3),
\end{split}\end{equation}
where
\begin{equation}\begin{split}
\delta_0=-\gamma_0-\gamma_1-\sum_{i=0}^N\frac{m_i}{\lambda_i-1},
\quad
\delta_1=\gamma_1+\sum_{i=0}^N\frac{m_i}{(1-\lambda_i)^2},
\quad
\delta_2=-\frac{\gamma_1}{2}
+\sum_{i=0}^N\frac{m_i(1+\lambda_i)}{2(1-\lambda_i)^3}.
\end{split}\end{equation}
In accordance with condition (\ref{b6}) for \( p = 2 \), we have \( \delta_0 = 0 \), \( \delta_1 = 1 \), and \( \delta_2 = 0 \), but this still leads to contradiction:
\begin{itemize}
\item If $N<0$, we would have \( \delta_1 = \gamma_1 = 1 \) but simultaneously \( \delta_2 = -\frac{\gamma_1}{2} = 0 \), which is a contradiction.
\item If $N \geq 0$, since $m_i>0$, $\lambda_i>1$ and $\gamma_1 \geq 0$, then $\delta_2=-\frac{\gamma_1}{2}
+\sum_{i=0}^N\frac{m_i(1+\lambda_i)}{2(1-\lambda_i)^3}<0$,  a contradiction.
\end{itemize}
The proof of the lemma is completed.
\end{proof}
\end{lemma}
\par
The above lemma indicates that
it is challenging to directly construct a high-order generating function \( w(\zeta) \)  that is also a Pick function and is analytic for $\zeta \in (-\infty,1)$ within the framework of CQ theory by assuming  \( w(\zeta) \)  is of the form in (\ref{b4}).  
We introduce a perturbation to the generating function \( w(\zeta) \) by assuming it has the following structure
\begin{equation}\begin{split}\label{b5}
w(\zeta)=\bigg(1+\bigg(\frac{1-\zeta}{\tau}\bigg)^\alpha G(\zeta)^{1-\alpha}\bigg)^{-\beta}.
\end{split}\end{equation}
This assumption is motivated by the fact that if \( -G(\zeta) \) is a Pick function, then \( w(\zeta) \) is also a Pick function.
\begin{lemma}
If $-G(\zeta)$ is a Pick function, then $w(\zeta)$ defined in (\ref{b5}) is a Pick function.
\begin{proof}
Let \(z = (1 - \zeta)/\tau\). For \(\zeta \in \mathbb{C}^+\), \(\text{Im}(\zeta) > 0\) implies \(\text{Im}(z) < 0\), so \(z\) lies strictly in the lower half-plane. Using the principal branch of the power function, \(z^\alpha\) is analytic and its argument satisfies \(\arg(z^\alpha) \in (-\alpha\pi, 0)\), hence \(z^\alpha\) remains strictly in the lower half-plane.

Given that \(-G(\zeta)\) is a Pick function with \(\text{Im}(-G(\zeta)) > 0\) (by nonconstancy), we have \(\text{Im}(G(\zeta)) < 0\). Then \(G(\zeta)^{1-\alpha}\) is analytic and its argument satisfies \(\arg(G(\zeta)^{1-\alpha}) \in (-(1-\alpha)\pi, 0)\), so it also lies strictly in the lower half-plane.

Now define \(H(\zeta) = z^\alpha G(\zeta)^{1-\alpha}\). Then \(H(\zeta)\) is analytic on \(\mathbb{C}^+\), and its argument is  
\[
\arg(H(\zeta)) = \arg(z^\alpha) + \arg(G(\zeta)^{1-\alpha}) \in (-\pi, 0),
\]  
so \(\text{Im}(H(\zeta)) < 0\). Thus, \(H(\zeta)\) is strictly in the lower half-plane.

Consider \(1 + H(\zeta)\). Since \(\text{Im}(1 + H(\zeta)) = \text{Im}(H(\zeta)) < 0\), it follows that \(1 + H(\zeta)\) is strictly in the lower half-plane and avoids the negative real axis (as its argument is in \((-\pi, 0)\)).

Finally,  
\[
w(\zeta) = (1 + H(\zeta))^{-\beta} = \exp\left( -\beta \log(1 + H(\zeta)) \right).
\]  
Because \(1 + H(\zeta)\) is in the lower half-plane, \(\arg(1 + H(\zeta)) \in (-\pi, 0)\), so  
\[
\arg(w(\zeta)) = -\beta \arg(1 + H(\zeta)) \in (0, \beta\pi) \subset (0, \pi),
\]  
which implies \(\text{Im}(w(\zeta)) > 0\). Moreover, \(w(\zeta)\) is analytic on \(\mathbb{C}^+\) as a composition of analytic functions that avoid branch cuts.
Therefore, \(w(\zeta)\) is a Pick function.
\end{proof}
\end{lemma}
\par
The following three steps result in a completely monotonic sequence  $(w_0, w_1, \cdots )$ such that the discrete convolution $\sum_{k=0}^n w_{n-k}(\tau)u^k$ approximates $\int_0^t \omega_{\alpha,\beta}(t-s)u(s)ds$ with second-order accuracy.
\\
\par
\textbf{Step I.}
We assume that $G(\zeta)$ is a rational function  as simply as possible, and that $-G(\zeta)$ is a Pick function.
Since we aim to constructing an approximation formula with second-order accuracy and  due to Lemma \ref{lem.3}, we set $G(\zeta)=-\gamma_1\zeta-\gamma_0$ where $\gamma_1 \geq 0$.
\par
\textbf{Step II.}
With $G(\zeta)$ defined above, in accordance with (\ref{b4}) and (\ref{b5}), we get $\delta(\zeta)=\tau^{-1}(1-\zeta)G(\zeta)^{\frac{1-\alpha}{\alpha}}$.
By requiring that $\delta(\zeta)$ satisfies the condition (\ref{b6}) with $p=2$, we can obtain the specific $G(\zeta)$.
\begin{lemma}\label{lem.4}
The function $\delta(\zeta)=\tau^{-1}(1-\zeta)(-\gamma_1\zeta-\gamma_0)^{\frac{1-\alpha}{\alpha}}$ satisfies the condition (\ref{b6}) with $p=2$ if and only if
$\gamma_0=-\frac{2-\alpha}{2-2\alpha}$, $\gamma_1=\frac{\alpha}{2-2\alpha}$.
\begin{proof}
We prove the lemma by matching the asymptotic expansion of  
\[
\frac{1}{\tau} \delta(e^{-\tau}) = \frac{1}{\tau}(1 - e^{-\tau})(-\gamma_1 e^{-\tau} - \gamma_0)^{\frac{1-\alpha}{\alpha}}
\]  
with the condition \( 1 + O(\tau^2) \) as \( \tau \to 0 \).
Let \( \tilde{\alpha} = \frac{1-\alpha}{\alpha} \). Expanding \( e^{-\tau} = 1 - \tau + \frac{\tau^2}{2} - \frac{\tau^3}{6} + O(\tau^4) \), we have  
\[
1 - e^{-\tau} = \tau - \frac{\tau^2}{2} + \frac{\tau^3}{6} + O(\tau^4),
\]  
\[
-\gamma_1 e^{-\tau} - \gamma_0 = -(\gamma_0 + \gamma_1) + \gamma_1 \tau - \frac{\gamma_1}{2} \tau^2 + \frac{\gamma_1}{6} \tau^3 + O(\tau^4).
\]
Let \( A = -(\gamma_0 + \gamma_1) \). The leading-order term is \( A^{\tilde{\alpha}} \), so we require  
\begin{equation}\begin{split}\label{b12}
A^{\tilde{\alpha}} = 1 \quad \Rightarrow \quad A = 1 \quad \Rightarrow \quad \gamma_0 + \gamma_1 = -1.
\end{split}\end{equation}
Under (\ref{b12}), we expand  
\[
(1 + \gamma_1 \tau - \frac{\gamma_1}{2} \tau^2 + \cdots)^{\tilde{\alpha}} 
= 1 + {\tilde{\alpha}} \gamma_1 \tau + \left[ -\frac{{\tilde{\alpha}} \gamma_1}{2} + \frac{{\tilde{\alpha}}({\tilde{\alpha}}-1)}{2} \gamma_1^2 \right] \tau^2 + O(\tau^3).
\]  
Multiplying by \( 1 - e^{-\tau} \) and dividing by \( \tau \), we obtain  
\[
\frac{1}{\tau} \delta(e^{-\tau}) = 1 + \left( {\tilde{\alpha}} \gamma_1 - \frac12 \right) \tau + O(\tau^2).
\]  
To eliminate the \( O(\tau) \) term, we require  
\[
{\tilde{\alpha}} \gamma_1 - \frac12 = 0 \quad \Rightarrow \quad \gamma_1 = \frac{1}{2{\tilde{\alpha}}} = \frac{\alpha}{2(1-\alpha)}. \tag{2}
\]
From (1) and (2),  
\[
\gamma_0 = -1 - \gamma_1 = -\frac{2 - \alpha}{2(1-\alpha)}.
\]

Substituting these values into the expansion confirms that the \( O(\tau) \) term vanishes and the \( O(\tau^2) \) term is bounded, yielding \( 1 + O(\tau^2) \). The converse follows from the uniqueness of asymptotic coefficients. 
\end{proof}
\end{lemma}
\par
\textbf{Step III.} With $\delta(\zeta)$ derived from \textbf{Step II}, verify that $\omega(\zeta)$ in (\ref{b5}) is analytic and nonnegative on $(-\infty,1)$, which further indicates $(w_0, w_1, \cdots )$ is completely monotonic by Lemma \ref{lem.2}.
\begin{lemma}\label{lem.5}
With $\delta(\zeta)$ in Lemma \ref{lem.4}, the generating function $w(\zeta)$ takes the form
\begin{equation}\begin{split}\label{b10}
w(\zeta)=
\bigg[
1+\bigg(\frac{1-\zeta}{\tau}\bigg)^\alpha
\bigg(\frac{2-\alpha}{2-2\alpha}\bigg)^{1-\alpha}
\bigg(1-\frac{\alpha}{2-\alpha}\zeta\bigg)^{1-\alpha}
\bigg]^{-\beta},\quad 0<\alpha,\beta<1,
\end{split}\end{equation}
and is analytic and nonnegative on $(-\infty,1)$.
\begin{proof}
For \(\zeta \in (-\infty, 1)\), we have \(1 - \zeta > 0\) and  
\[
1 - \frac{\alpha}{2 - \alpha} \zeta > 1 - \frac{\alpha}{2 - \alpha} = \frac{2 - 2\alpha}{2 - \alpha} > 0.
\]  
Thus, all factors inside the brackets are positive, and the expression  
\[
B(\zeta) := \left( \frac{1 - \zeta}{\tau} \right)^\alpha \left( \frac{2 - \alpha}{2 - 2\alpha} \right)^{1 - \alpha} \left( 1 - \frac{\alpha}{2 - \alpha} \zeta \right)^{1 - \alpha}
\]  
is a composition of real-analytic functions on \((-\infty, 1)\), hence real-analytic there. Since \(B(\zeta) > 0\), we have \(1 + B(\zeta) > 1\), and \(w(\zeta) = [1 + B(\zeta)]^{-\beta}\) is also real-analytic (as composition of real-analytic functions) and strictly positive.  
\end{proof}
\end{lemma}
\subsection{Energy-decay schemes for the H-N model}
Introduce the symbol $\partial_\tau u^{n-\frac{1}{2}}=\frac{u^{n}-u^{n-1}}{\tau}$ and $u^{n-\frac{1}{2}}=\frac{1}{2}(u^n+u^{n-1})$.
The time semi-discrete scheme reads that finding \( \boldsymbol E^n \in H_0(\mathrm{curl}, \Omega) \)  and  \( \boldsymbol H^n, \boldsymbol P^n \in \boldsymbol L^2(\Omega) \)  fulfilling
\begin{align}
    &\epsilon_\infty (\partial_\tau \boldsymbol E^{n-\frac{1}{2}}, \boldsymbol\phi) + (\partial_\tau \boldsymbol P^{n-\frac{1}{2}}, \boldsymbol\phi) - (\boldsymbol H^{n-\frac{1}{2}}, \nabla \times \boldsymbol\phi) = 0, 
    & \forall \boldsymbol\phi \in H_0(\mathrm{curl}, \Omega), \tag{3.1}\label{b7} \\
    &(\partial_\tau \boldsymbol H^{n-\frac{1}{2}}, \boldsymbol\psi) + (\nabla \times \boldsymbol E^{n-\frac{1}{2}}, \boldsymbol\psi) = 0, 
    & \forall \boldsymbol\psi \in \boldsymbol L^2(\Omega), \tag{3.2}\label{b8}  \\
    &(\boldsymbol P^n,\boldsymbol\varphi)- \Delta \epsilon \sum_{k=0}^n w_{n-k} (\boldsymbol E^k,\boldsymbol\varphi)=0, 
    &  \forall \boldsymbol\varphi \in \boldsymbol L^2(\Omega) \tag{3.3}\label{b9} 
\end{align}
where the sequence $(w_0, w_1, \cdots )$ is generated from (\ref{b10}).
\begin{lemma}\label{lem.11}
If $(w_0, w_1, \cdots )$ is a completely monotonic sequence, there holds
\begin{equation}\begin{split}
(\boldsymbol E^{n-\frac{1}{2}},\boldsymbol P^n-\boldsymbol P^{n-1})
\geq \frac{\Delta_\epsilon}{2}
\bigg(
\sum_{k=0}^n w_{n-k}\|\boldsymbol E^k\|^2
-
\sum_{k=0}^{n-1} w_{n-1-k}\|\boldsymbol E^k\|^2
\bigg).
\end{split}\end{equation}
\begin{proof}
Assume $w_i=0$ for $i<0$.
Using (\ref{b9}), we have
\begin{equation}\begin{split}
(\boldsymbol E^{n-\frac{1}{2}},\boldsymbol P^n-\boldsymbol P^{n-1})=
\Delta_\epsilon
\sum_{k=0}^n(w_{n-k}-w_{n-1-k})
(\boldsymbol E^k, \boldsymbol E^{n-\frac{1}{2}}).
\end{split}\end{equation}
Since the sequence $(w_0, w_1, \cdots )$ is completely monotonic, $w_i\geq 0$ for $i \geq 0$ and $w_{i_1} \geq w_{i_2}$ for $i_1<i_2$.
Based on the Cauchy-Schwarz inequality and the Young inequality, we get
\begin{equation}\begin{split}
(\boldsymbol E^{n-\frac{1}{2}},\boldsymbol P^n-\boldsymbol P^{n-1})
&\geq
\frac{\Delta_\epsilon}{2}
\sum_{k=0}^n(w_{n-k}-w_{n-1-k})
(\|\boldsymbol E^k\|^2+\| \boldsymbol E^{n-\frac{1}{2}}\|^2)
\\&=
\frac{\Delta_\epsilon}{2}w_n \| \boldsymbol E^{n-\frac{1}{2}}\|^2
+ \frac{\Delta_\epsilon}{2}
\bigg(
\sum_{k=0}^n w_{n-k}\|\boldsymbol E^k\|^2
-
\sum_{k=0}^{n-1} w_{n-1-k}\|\boldsymbol E^k\|^2
\bigg),
\end{split}\end{equation}
which completes the proof of the lemma.
\end{proof}
\end{lemma}
\begin{theorem}
If the sequence $(w_0, w_1, \cdots )$ is completely monotonic, then the numerical scheme (\ref{b7})-(\ref{b9}) preserving  the discrete energy-decay property $\mathbb{E}^n \leq \mathbb{E}^{n-1}, n \geq 1$ where
\begin{equation}\begin{split}
\mathbb{E}^n=
\epsilon_\infty \|\boldsymbol{E}^n\|^2+\|\boldsymbol{H}^n\|^2 
+\Delta_\epsilon\sum_{k=0}^n w_{n-k}\|\boldsymbol E^k\|^2.
\end{split}\end{equation}
\begin{proof}
By setting $\boldsymbol\phi = \boldsymbol E^{n-\frac{1}{2}}$ in equation (\ref{b7}) and $\boldsymbol\psi = \boldsymbol H^{n-\frac{1}{2}}$ in equation (\ref{b8}), and adding the two resulting equations, we obtain
\begin{equation}\begin{split}
\epsilon_\infty (\partial_\tau \boldsymbol E^{n-\frac{1}{2}}, \boldsymbol E^{n-\frac{1}{2}}) 
+ (\partial_\tau \boldsymbol P^{n-\frac{1}{2}}, \boldsymbol E^{n-\frac{1}{2}})
+(\partial_\tau \boldsymbol H^{n-\frac{1}{2}}, \boldsymbol H^{n-\frac{1}{2}})
 = 0.
\end{split}\end{equation}
Combining Lemma \ref{lem.11} and the above equation, one can get readily that
\[
\mathbb{E}^n \leq \mathbb{E}^{n-1}, n \geq 1.
\]
\end{proof}
\end{theorem}
\begin{remark}\label{rem.1}
The main objective of this paper is to construct high-order completely monotone sequences that approximate the completely monotone kernel function $\omega_{\alpha,\beta}(t)$, thereby deriving numerical schemes for the H-N model that preserve the energy decay property. For spatial discretization, we adopt \textit{Nédélec} elements directly and omit the convergence analysis of the fully discrete scheme, as such analytical techniques are well-established and can be found, for example, in \cite{baoli2023discrete}.
Although the original model (1.1)-(1.5) is formulated in three dimensions, we adopt a two-dimensional Maxwell framework for computational convenience. In this 2D setting, the curl operators are defined as follows:
\[
\nabla \times H = \left( \frac{\partial H}{\partial y}, -\frac{\partial H}{\partial x} \right)^\mathrm{T}, \quad \nabla \times \boldsymbol{E} = \frac{\partial E_2}{\partial x} - \frac{\partial E_1}{\partial y}.
\]
Let \( \mathcal{T}_h \) be a family of regular rectangular partitions of the domain \( \Omega = (0,1)^2 \) with maximum mesh size \( h \). For any integer \( k \geq 1 \), we define the following mixed finite element spaces \cite{li2012time}:

\[
\begin{aligned}
V_h &= \left\{ \psi_h \in L^2(\Omega) : \psi_h|_K \in Q_{k-1,k-1}, ~ \forall K \in \mathcal{T}_h \right\}, \\
\boldsymbol{U}_h &= \left\{ \boldsymbol{\phi}_h = (\phi_{h1}, \phi_{h2})^\mathrm{T} \in H(\mathrm{curl};\Omega) : 
\phi_{h1}|_K \in Q_{k-1,k}, ~ \phi_{h2}|_K \in Q_{k,k-1}, ~ \forall K \in \mathcal{T}_h \right\}, \\
\boldsymbol{U}_h^0 &= \left\{ \boldsymbol{\phi}_h \in \boldsymbol{U}_h : 
\phi_{h1}|_{y=0} = \phi_{h1}|_{y=1} = 0, ~ \phi_{h2}|_{x=0} = \phi_{h2}|_{x=1} = 0 \right\}.
\end{aligned}
\]
where
\begin{equation}\label{I.6}\begin{split}
Q_{i,j}&=\{\text{polynomials of maximum degrees $i$, $j$ in $x, y$, respectively}\}.
\end{split}
\end{equation}
Note that for any \( \boldsymbol{u}_h \in \boldsymbol{U}_h^0 \), the boundary condition \( \boldsymbol{n} \times \boldsymbol{u}_h = \boldsymbol{0} \) is satisfied on \( \partial \Omega \).
The fully discrete approximation of the system (\ref{a1})-(\ref{a3}) is then formulated as follows: at each time step, find \( \boldsymbol{E}_h^n, \boldsymbol{P}_h^n \in \boldsymbol{U}_h^0 \) and \( H_h^n \in V_h \) such that

\[
\begin{aligned}
\epsilon_\infty \left( \partial_\tau^{n-\frac{1}{2}} \boldsymbol{E}_h, \boldsymbol{\phi}_h \right)
+ \left( \partial_\tau^{n-\frac{1}{2}} \boldsymbol{P}_h, \boldsymbol{\phi}_h \right)
- \left( {H}_h^{n-\frac{1}{2}}, \nabla \times \boldsymbol{\phi}_h \right) &= 0,
\quad \forall \boldsymbol{\phi}_h \in \boldsymbol{U}_h^0, \\
 \left( \partial_\tau^{n-\frac{1}{2}} H_h, \psi_h \right)
+ \left( \nabla \times {\boldsymbol{E}}_h^{n-\frac{1}{2}}, \psi_h \right) &= 0,
\quad \forall \psi_h \in V_h, \\
\left(  \boldsymbol{P}_h^n, \boldsymbol{\varphi}_h \right)
-
\Delta_\epsilon
\sum_{k=0}^n w_{n-k}
 \left({\boldsymbol{E}}_h^{k}, \boldsymbol{\varphi}_h \right) &= 0,
\quad \forall \boldsymbol{\psi}_h \in \boldsymbol{U}_h.
\end{aligned}
\]

\end{remark}

\section{Numerical tests}
In this section, we perform numerical experiments to validate the theoretical predictions for the completely monotonic sequence, energy-decay property, and convergence rates.
\subsection{Validation of the completely monotonic property}
Following (\ref{b00}), we define, for a sequence $(w_0, w_1, \cdots )$, the index
\[
\text{Index}_k (\alpha,\beta,\tau)=\min_{0\leq j \leq J}(I-S)^k w_j(\alpha,\beta,\tau),\quad k=0,1,\cdots.
\]
Then, for given $\alpha,\beta,\tau$, the sequence is completely monotonic if and only if $\text{Index}_k(\alpha,\beta,\tau) \geq 0$ for all $k \geq 0$ (with $J=\infty$). Specifically, $\text{Index}_0 \geq 0$ implies $w_j \geq 0$, $\text{Index}_1 \geq 0$ implies $w_{j+1} - w_j \leq 0$, $\text{Index}_2 \geq 0$ implies $w_{j+2} - 2w_{j+1} + w_j \geq 0$, and so forth. These conditions constitute a discrete analogue of the continuous property (\ref{b000}).
Introduce the function
\[
\rho(x)=
\begin{cases}
1,&\quad x\geq 0,
\\
0,&\quad x<0.
\end{cases}
\]
\par
\begin{figure}[htbp]
\centering
\subfigure[]{
\begin{minipage}[t]{0.5\linewidth}
\centering
\includegraphics[width=1\textwidth]{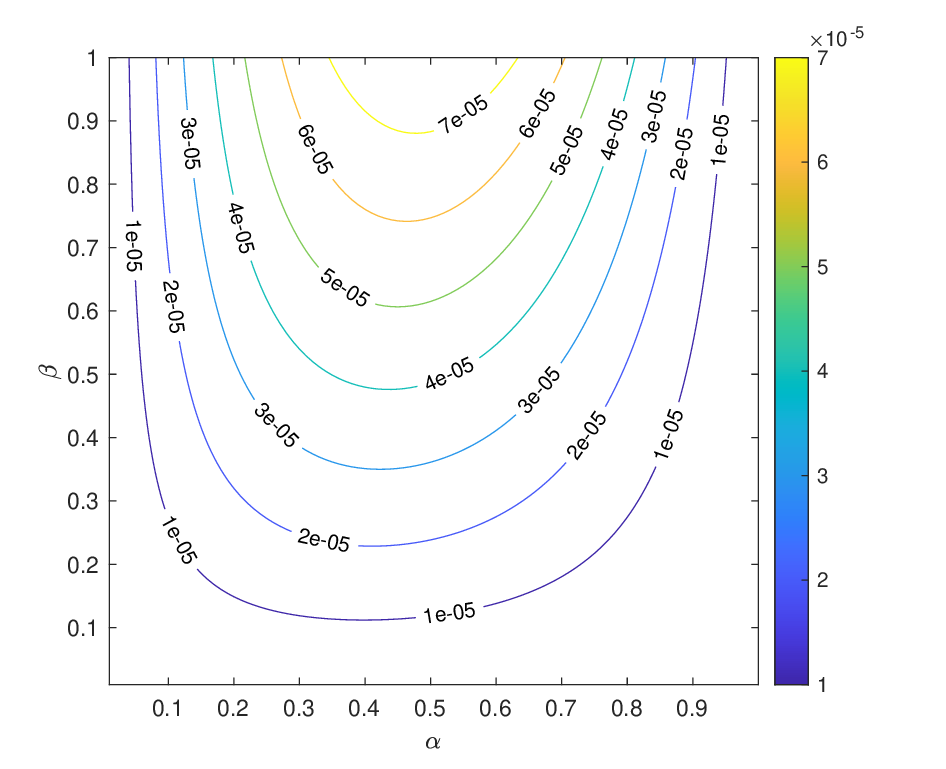}
\end{minipage}%
}%
\subfigure[]{
\begin{minipage}[t]{0.5\linewidth}
\centering
\includegraphics[width=1\textwidth]{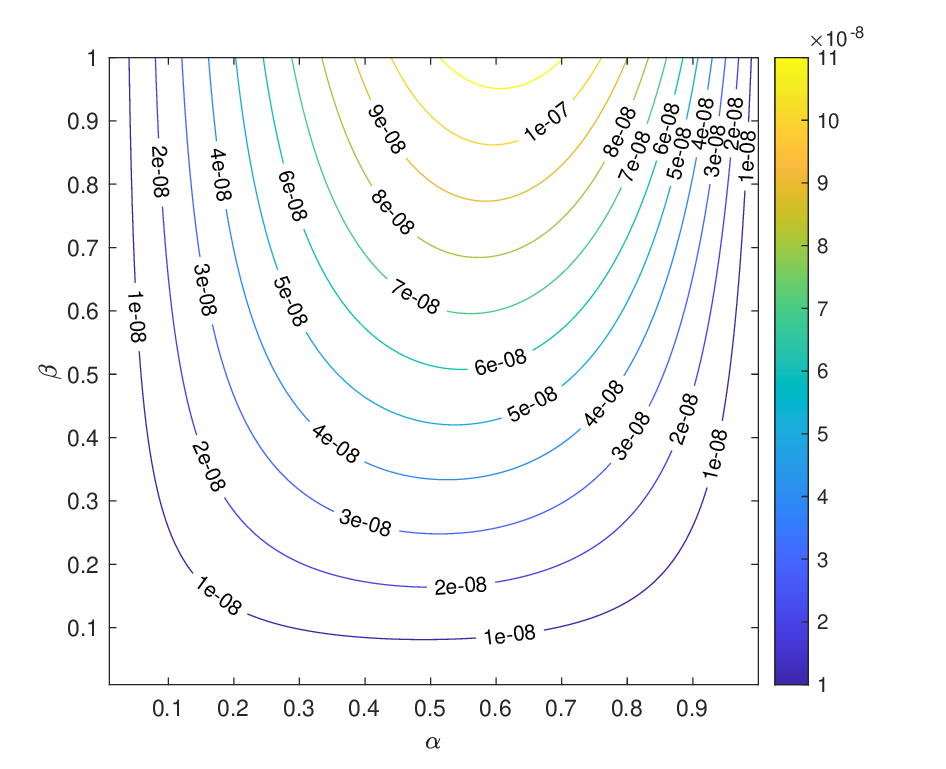}
\end{minipage}%
}%
\\
\subfigure[]{
\begin{minipage}[t]{0.5\linewidth}
\centering
\includegraphics[width=1\textwidth]{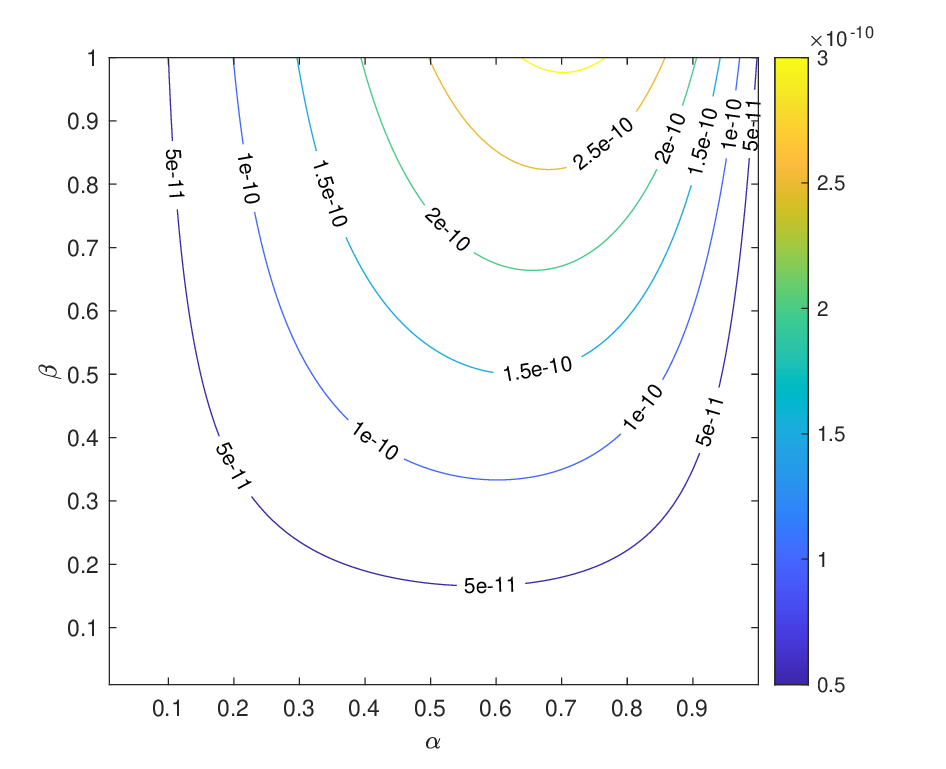}
\end{minipage}%
}%
\subfigure[]{
\begin{minipage}[t]{0.5\linewidth}
\centering
\includegraphics[width=1\textwidth]{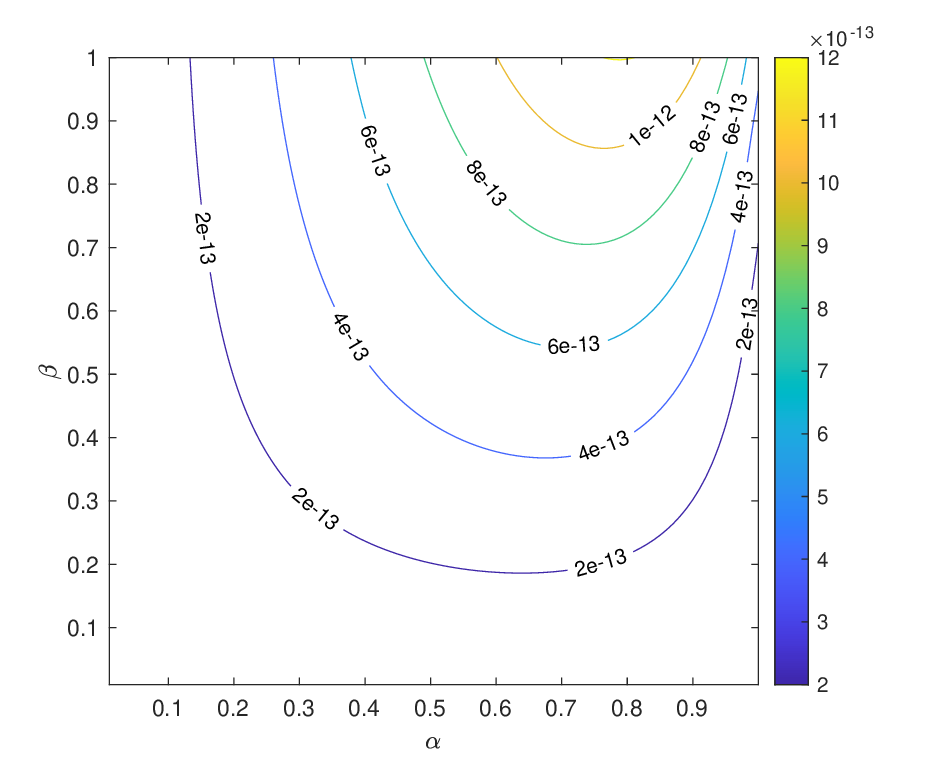}
\end{minipage}%
}%
\centering
\caption{Validation of the completely monotonic property for the sequence $(w_0, w_1, \cdots )$ with $J=1000, \tau=0.01$: (a) $\text{Index}_0$; (b) $\text{Index}_1$; (c) $\text{Index}_2$; (d) $\text{Index}_3$.}\label{C1}
\end{figure}

\begin{figure}[htbp]
\centering
\subfigure[]{
\begin{minipage}[t]{0.5\linewidth}
\centering
\includegraphics[width=1\textwidth]{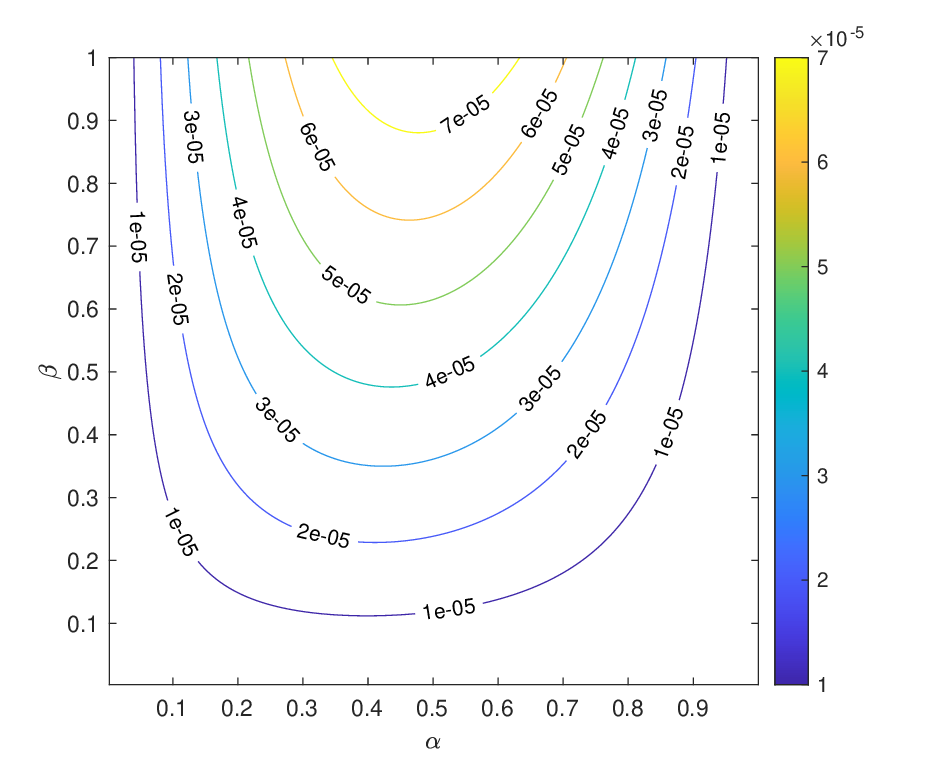}
\end{minipage}%
}%
\subfigure[]{
\begin{minipage}[t]{0.5\linewidth}
\centering
\includegraphics[width=1\textwidth]{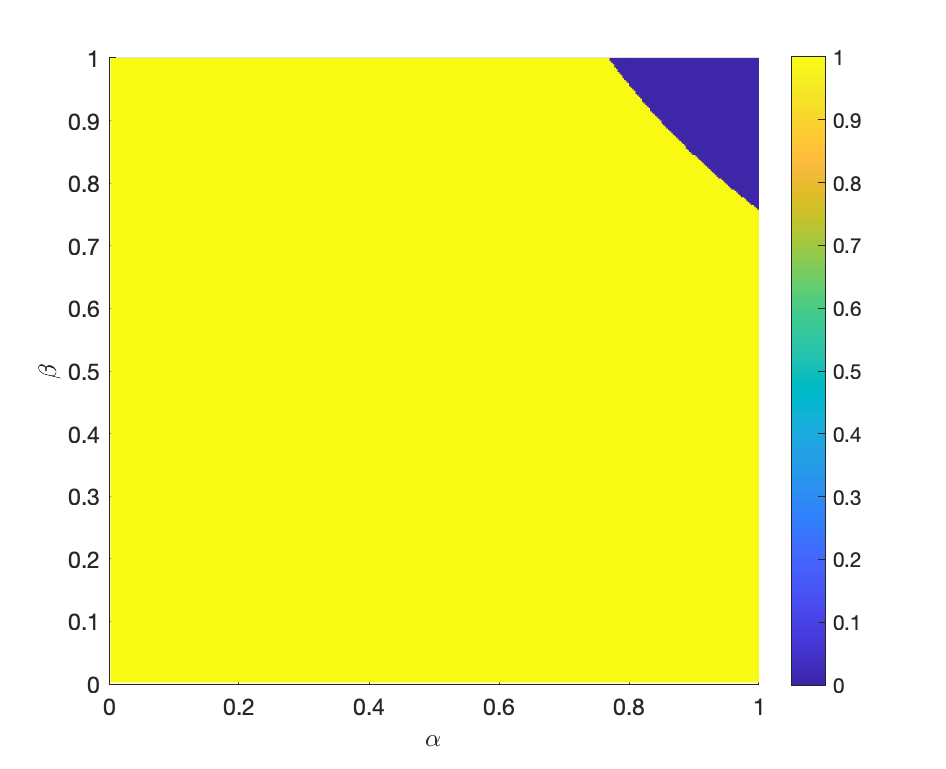}
\end{minipage}%
}%
\\
\subfigure[]{
\begin{minipage}[t]{0.5\linewidth}
\centering
\includegraphics[width=1\textwidth]{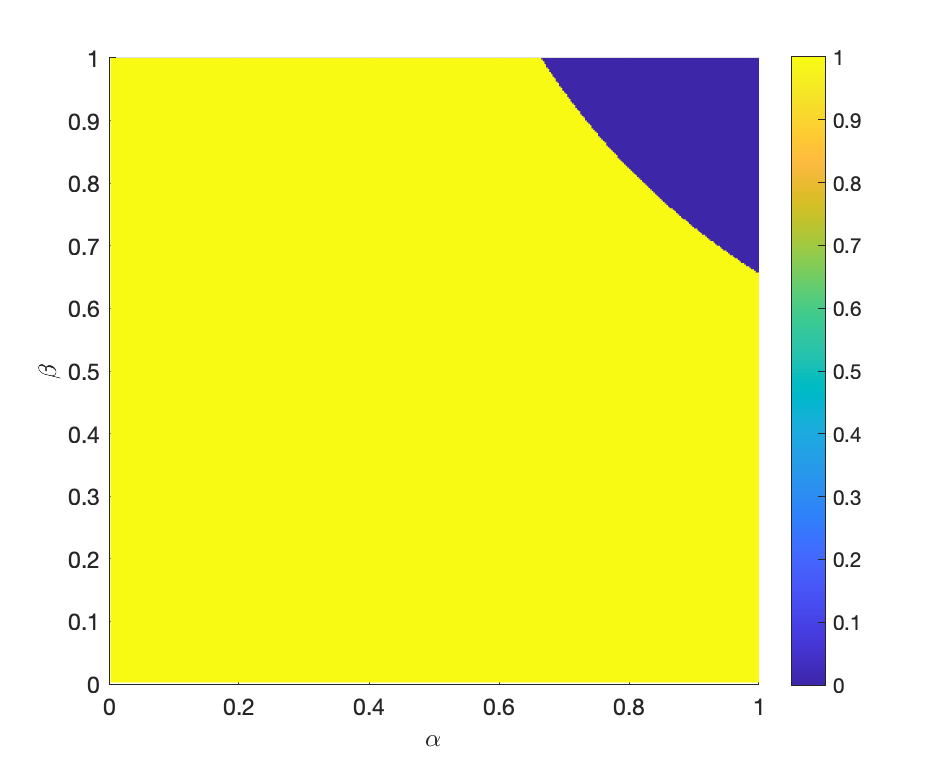}
\end{minipage}%
}%
\subfigure[]{
\begin{minipage}[t]{0.5\linewidth}
\centering
\includegraphics[width=1\textwidth]{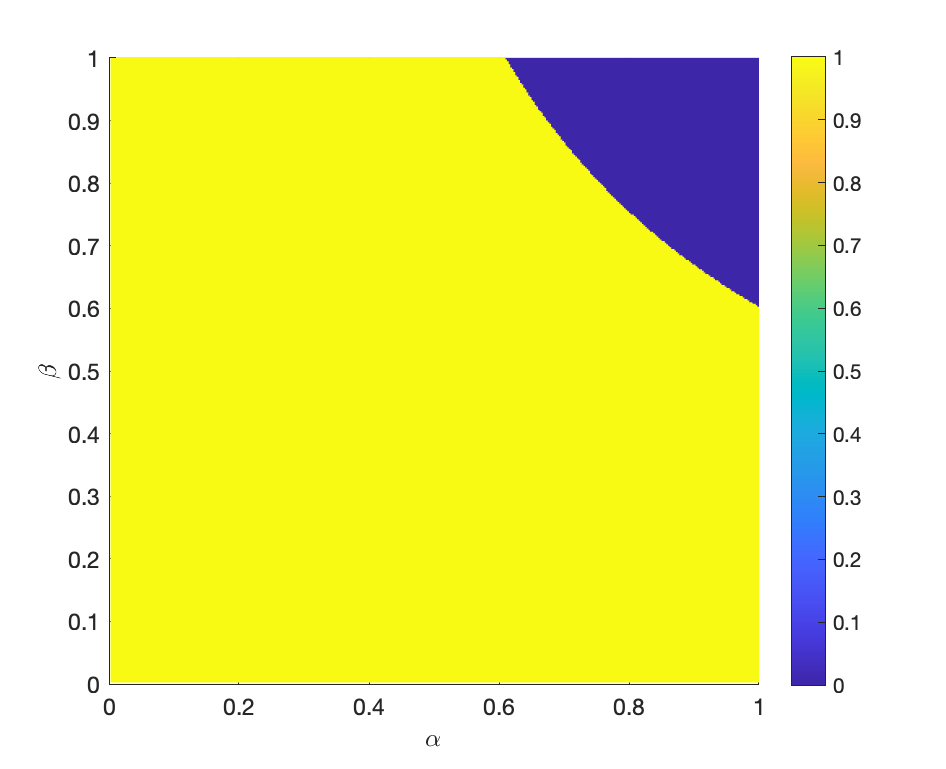}
\end{minipage}%
}%
\centering
\caption{Illustration of the non-completely monotonic property for the sequence $(w_0, w_1, \cdots )$ generated by the BDF-2 with $J=1000, \tau=0.01$: (a) the value of $\text{Index}_0$; (b) the value of $\rho(\text{Index}_1)$; (c) the value of $\rho(\text{Index}_2)$; (d) the value of $\rho(\text{Index}_3)$.}\label{C2}
\end{figure}

 
In Fig. \ref{C1}, the non-negativity of the sequence generated from (\ref{b10}) is verified by plotting $\text{Index}_k$ for $k=0,1,2,3$ with $J=1000$ and $\tau=0.01$.
It is evident that all examined terms in the sequence satisfy the non-negativity condition, as illustrated.
For comparison, Fig. \ref{C2} displays $\text{Index}_k$ for the sequence generated from (\ref{b4}) with $\delta(\zeta)$ taken as the generating function of the BDF-2 method.
Although $\text{Index}_0 \geq 0$ in this case, $\text{Index}_k$ (for $k=1,2,3$) clearly becomes negative for certain values of $\alpha$ and $\beta$. Moreover, as $k$ increases, an expanding set of parameter pairs $(\alpha,\beta)$ results in negative $\text{Index}_k$.
These results demonstrate that the BDF-2 scheme cannot generate a completely monotonic sequence.

\subsection{H-N Model with Source Terms}
We consider the spatial domain $\Omega = (0,1)^2$ and the final time $T = 1$. 
The physical parameters are chosen as  
\[
\epsilon_{\infty} = \Delta \epsilon = 1.
\]  
To quantitatively examine the convergence behavior, we introduce the following set of smooth analytical solutions:  
\begin{equation}\label{N.1}
\begin{aligned}
\boldsymbol{E}(x,y,t) &= t^{3}
\begin{pmatrix}
\displaystyle(x^2 + 1) \sin(\pi y) \\
\displaystyle\sin(\pi x)\left(y - \frac{1}{2}\right)
\end{pmatrix}, \\
\boldsymbol{P}(x,y,t) &= (1 - e^{-t})
\begin{pmatrix}
\displaystyle(x^2 + 1) y(y - 1) \\
\displaystyle x(x - 1)\left(y - \frac{1}{2}\right)
\end{pmatrix}, \\
H(x,y,t) &= \displaystyle e^{-t}(x^3 + 1)(y^3 + 1).
\end{aligned}
\end{equation}  
Using \eqref{b1.2} with the prescribed solutions, one can analytically derive the corresponding source terms.
Adopting the mixed finite element spaces in Remark \ref{rem.1} with $k=2$ and employing the error estimation methodology for the fully discrete scheme from \cite{baoli2023discrete}, we derive the following error estimates for the exact solution:
\begin{equation}\begin{split}
\max_n\|\boldsymbol E^n-\boldsymbol E_h^n\|
+\max_n\| H^n- H_h^n\|
+\max_n\|\boldsymbol P^n-\boldsymbol P_h^n\|
=O(\tau^\ell+h^2),
\end{split}\end{equation}
where $\ell=1$ if $\alpha=1$, and 2 otherwise.
For notational simplicity, we define $\text{Error}(\boldsymbol{E}) = \max_{0 \le n \le N} \| \boldsymbol{E}(t_n) - \boldsymbol{E}_h^n \|$, with the errors for $\boldsymbol{P}$ and $H$ defined similarly.
\par
To validate the theoretical error estimates, numerical experiments on the temporal convergence rate were performed with a fixed spatial mesh size \(h = \frac{\sqrt{2}}{100}\).
The results are presented in Table \ref{tab1} (for the H-N model) and Table \ref{tab2} (for the C-C model with \(\beta = 1\)). 
For the H-N model (Table \ref{tab1}), under a wide range of parameter pairs \((\alpha, \beta)\), most cases also exhibit convergence rates close to second order. Only when the parameters approach 1 (e.g., \((0.9,0.9)\)) do some rates show a slight reduction; however, as the time step is further refined, the rates stabilize and remain consistent with the theoretical predictions. 
For the C-C model (Table \ref{tab2}), the temporal convergence rates of the errors for \(\boldsymbol E\), \(H\), and \(\boldsymbol P\) remain consistently around 2.00 for various values of \(\alpha\), confirming clear second‑order accuracy. 
Overall, the numerical results align well with the theoretical error estimates, demonstrating that the proposed fully discrete scheme achieves second‑order temporal convergence for a broad span of fractional parameters \(\alpha\) and \(\beta\).

\begin{table}[]
\centering
\caption{Temporal error and convergence rates for the H-N model with $h=\frac{\sqrt{2}}{100}$.}\label{tab1}
{\renewcommand{\arraystretch}{1.2}
\begin{tabular}{cccccccc}
\hline
$(\alpha,\beta)$           & $\tau$ & $\text{Error}(\boldsymbol{E})$ & Rates & $\text{Error}({H})$ & Rates & $\text{Error}(\boldsymbol{P})$ & Rates \\ \hline
\multirow{3}{*}{(0.1,0.1)} & 1/5    & 4.4248E-03                     &       & 9.9517E-03          &       & 4.2998E-03                     &       \\
                           & 1/10   & 1.1004E-03                     & 2.01  & 2.5133E-03          & 1.99  & 1.0724E-03                     & 2.00  \\
                           & 1/20   & 2.7539E-04                     & 2.00  & 6.2990E-04          & 2.00  & 2.6807E-04                     & 2.00  \\ \hline
\multirow{3}{*}{(0.1,0.5)} & 1/5    & 4.4878E-03                     &       & 1.0154E-02          &       & 3.8145E-03                     &       \\
                           & 1/10   & 1.1275E-03                     & 1.99  & 2.5581E-03          & 1.99  & 9.5659E-04                     & 2.00  \\
                           & 1/20   & 2.8143E-04                     & 2.00  & 6.4051E-04          & 2.00  & 2.3996E-04                     & 2.00  \\ \hline
\multirow{3}{*}{(0.1,0.9)} & 1/5    & 4.6493E-03                     &       & 1.0369E-02          &       & 3.3134E-03                     &       \\
                           & 1/10   & 1.1636E-03                     & 2.00  & 2.6073E-03          & 1.99  & 8.3454E-04                     & 1.99  \\
                           & 1/20   & 2.9092E-04                     & 2.00  & 6.5243E-04          & 2.00  & 2.0990E-04                     & 1.99  \\ \hline
\multirow{3}{*}{(0.5,0.1)} & 1/5    & 3.9818E-03                     &       & 9.2314E-03          &       & 5.1533E-03                     &       \\
                           & 1/10   & 9.8514E-04                     & 2.02  & 2.3087E-03          & 2.00  & 1.3185E-03                     & 1.97  \\
                           & 1/20   & 2.4655E-04                     & 2.00  & 5.7630E-04          & 2.00  & 3.3355E-04                     & 1.98  \\ \hline
\multirow{3}{*}{(0.5,0.5)} & 1/5    & 2.9901E-03                     &       & 7.4346E-03          &       & 7.9938E-03                     &       \\
                           & 1/10   & 6.9024E-04                     & 2.12  & 1.7865E-03          & 2.06  & 2.1246E-03                     & 1.91  \\
                           & 1/20   & 1.6872E-04                     & 2.03  & 4.3876E-04          & 2.03  & 5.4904E-04                     & 1.95  \\ \hline
\multirow{3}{*}{(0.5,0.9)} & 1/5    & 2.9449E-03                     &       & 6.7175E-03          &       & 9.1655E-03                     &       \\
                           & 1/10   & 6.3150E-04                     & 2.22  & 1.5667E-03          & 2.10  & 2.4692E-03                     & 1.89  \\
                           & 1/20   & 1.4654E-04                     & 2.11  & 3.8012E-04          & 2.04  & 6.4186E-04                     & 1.94  \\ \hline
\multirow{3}{*}{(0.9,0.1)} & 1/10   & 4.3776E-04                     &       & 1.2823E-03          &       & 3.0844E-03                     &       \\
                           & 1/20   & 6.6535E-05                     & 2.72  & 2.5881E-04          & 2.31  & 9.0264E-04                     & 1.77  \\
                           & 1/40   & 1.3910E-05                     & 2.26  & 5.8553E-05          & 2.14  & 2.4625E-04                     & 1.87  \\ \hline
\multirow{3}{*}{(0.9,0.5)} & 1/10   & 2.3742E-03                     &       & 2.6582E-03          &       & 8.5523E-03                     &       \\
                           & 1/20   & 6.7521E-04                     & 1.81  & 9.9154E-04          & 1.42  & 2.6235E-03                     & 1.70  \\
                           & 1/40   & 1.7778E-04                     & 1.93  & 2.9652E-04          & 1.74  & 7.3192E-04                     & 1.84  \\ \hline
\multirow{3}{*}{(0.9,0.9)} & 1/10   & 3.8540E-03                     &       & 3.7199E-03          &       & 1.0390E-02                     &       \\
                           & 1/20   & 1.1060E-03                     & 1.80  & 1.3786E-03          & 1.43  & 3.2377E-03                     & 1.68  \\
                           & 1/40   & 2.8746E-04                     & 1.94  & 4.1593E-04          & 1.73  & 9.1175E-04                     & 1.83  \\ \hline
\end{tabular}
}
\end{table}

\begin{table}[]
\centering
\caption{Temporal error and convergence rates for the C-C model ($\beta=1$) with $h=\frac{\sqrt{2}}{100}$.}\label{tab2}
{\renewcommand{\arraystretch}{1.2}
\begin{tabular}{ccllllll}
\hline
$\alpha$             & $\tau$ & \multicolumn{1}{c}{$\text{Error}(\boldsymbol{E})$} & \multicolumn{1}{c}{Rates} & \multicolumn{1}{c}{$\text{Error}({H})$} & \multicolumn{1}{c}{Rates} & \multicolumn{1}{c}{$\text{Error}(\boldsymbol{P})$} & \multicolumn{1}{c}{Rates} \\ \hline
\multirow{3}{*}{0.1} & 1/10   & 1.1729E-03                                         &                           & 2.6196E-03                              &                           & 8.0415E-04                                         &                           \\
                     & 1/20   & 2.9376E-04                                         & 2.00                      & 6.5544E-04                              & 2.00                      & 2.0237E-04                                         & 1.99                      \\
                     & 1/40   & 7.3408E-05                                         & 2.00                      & 1.6406E-04                              & 2.00                      & 5.0808E-05                                         & 1.99                      \\ \hline
\multirow{3}{*}{0.5} & 1/10   & 6.4840E-04                                         &                           & 1.5446E-03                              &                           & 2.4994E-03                                         &                           \\
                     & 1/20   & 1.4881E-04                                         & 2.12                      & 3.7391E-04                              & 2.05                      & 6.5029E-04                                         & 1.94                      \\
                     & 1/40   & 3.5851E-05                                         & 2.05                      & 9.2578E-05                              & 2.01                      & 1.6590E-04                                         & 1.97                      \\ \hline
\multirow{3}{*}{0.9} & 1/10   & 4.0811E-03                                         &                           & 3.7257E-03                              &                           & 1.0454E-02                                         &                           \\
                     & 1/20   & 1.1863E-03                                         & 1.78                      & 1.3927E-03                              & 1.42                      & 3.2717E-03                                         & 1.68                      \\
                     & 1/40   & 3.0751E-04                                         & 1.95                      & 4.2311E-04                              & 1.72                      & 9.2405E-04                                         & 1.82                      \\ \hline
\end{tabular}
}
\end{table}

\subsection{H-N Model with Zero Source Terms}
We validate the discrete energy-decay property of our presented scheme for the H-N model (\ref{a1})-(\ref{a3}) in this part.
We adopt the following initial condition
\begin{equation}\label{N.2}
\begin{aligned}
\boldsymbol{E}_0(x,y) &=
\begin{pmatrix}
\displaystyle(x^2 + 1) \sin(\pi y) \\
\displaystyle\sin(\pi x)\left(y - \frac{1}{2}\right)
\end{pmatrix}, \\
\boldsymbol{P}_0(x,y) &= 
\boldsymbol 0, \\
H_0(x,y) &= (x^3 + 1)(y^3 + 1).
\end{aligned}
\end{equation}  
The system is solved with zero source terms over the spatial domain \(\Omega = (0,1)^2\) up to a suitably chosen final time $T=1$. To examine the influence of the fractional parameters \(\alpha\) and \(\beta\) on the energy dissipation, we perform a series of numerical tests with \(\beta\) fixed at four representative values: \(\beta = 0.1, 0.4, 0.7, 1\). For each fixed \(\beta\), the fractional order \(\alpha\) is varied within the set \(\{0.1, 0.3, 0.5, 0.7, 0.9\}\).  
\par
\begin{figure}[htbp]
\centering
\subfigure[]{
\begin{minipage}[t]{0.5\linewidth}
\centering
\includegraphics[width=1\textwidth]{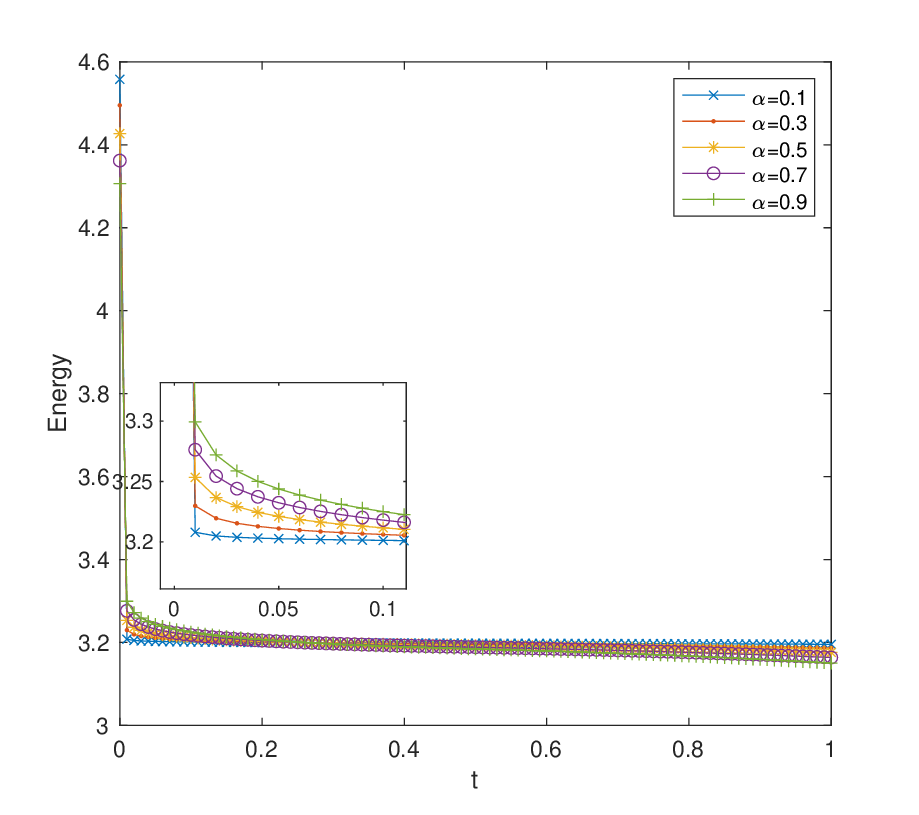}
\end{minipage}%
}%
\subfigure[]{
\begin{minipage}[t]{0.5\linewidth}
\centering
\includegraphics[width=1\textwidth]{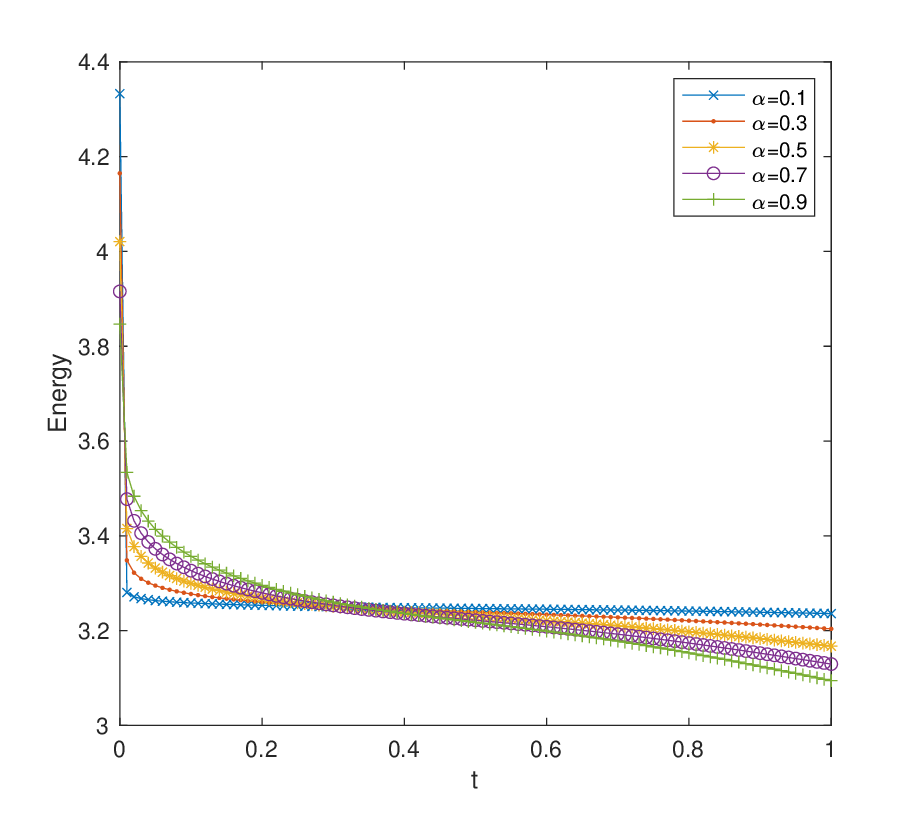}
\end{minipage}%
}%
\\
\subfigure[]{
\begin{minipage}[t]{0.5\linewidth}
\centering
\includegraphics[width=1\textwidth]{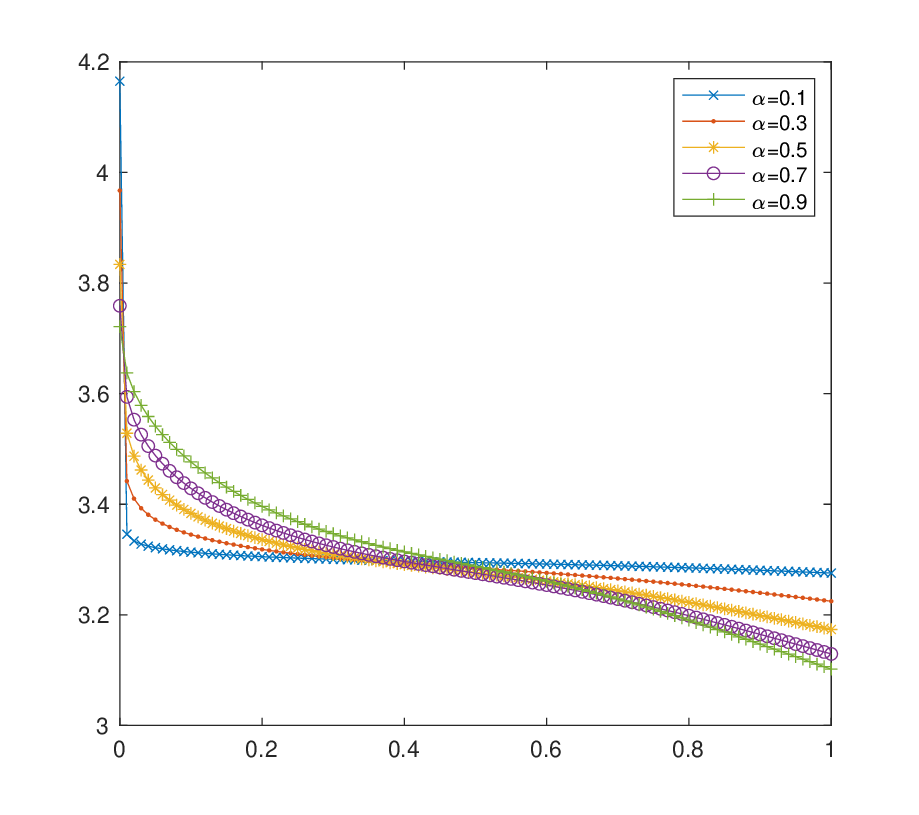}
\end{minipage}%
}%
\subfigure[]{
\begin{minipage}[t]{0.5\linewidth}
\centering
\includegraphics[width=1\textwidth]{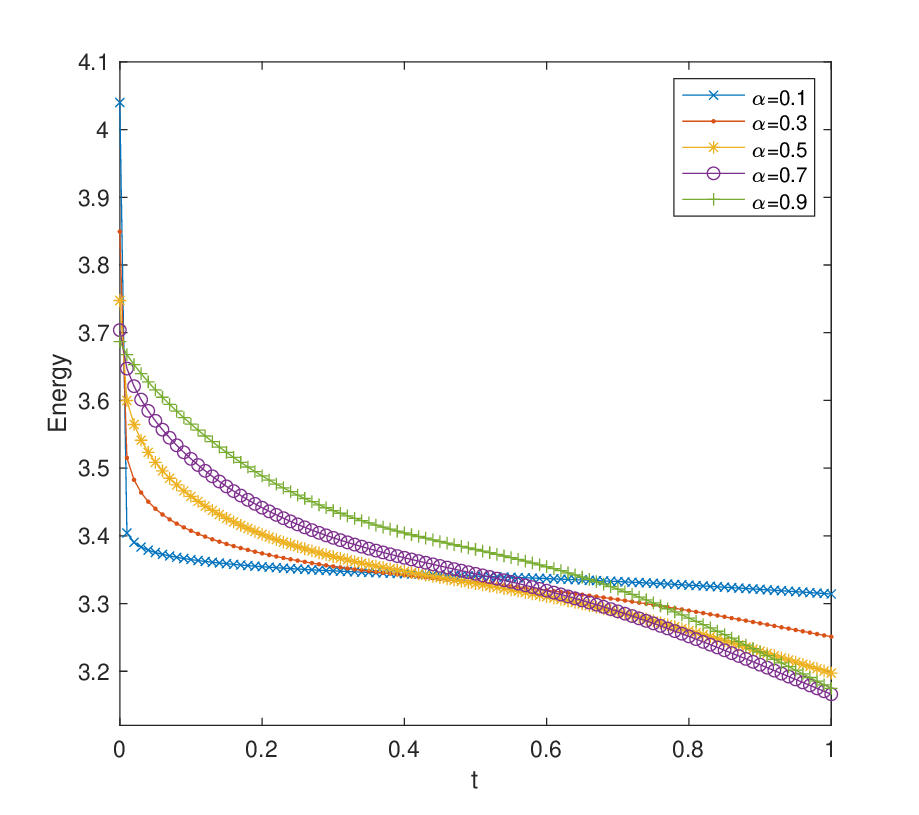}
\end{minipage}%
}%
\centering
\caption{Validation of energy decay characteristics across varying $\alpha$ and $\beta$: (a) $\beta=0.1$; (b) $\beta=0.4$; (c) $\beta=0.7$; (d) $\beta=1$. For each $\beta$, the parameter $\alpha$ is selected from the set $\alpha = 0.1, 0.3, 0.5, 0.7, 0.9$.}\label{C3}
\end{figure}

The evolution of the discrete energy (defined in accordance with the scheme) is monitored during the simulation. Fig. \ref{C3} displays the energy decay curves for the different parameter combinations. In each subfigure, corresponding to one value of \(\beta\), the curves for various \(\alpha\) are plotted together. The results clearly illustrate that the discrete energy decreases monotonically in time for all tested pairs \((\alpha, \beta)\), confirming that the proposed scheme preserves the energy-decay property of the continuous H-N model. Moreover, the rate of energy decay is observed to depend on both fractional orders: generally, a smaller \(\alpha\) or a smaller \(\beta\) leads to a faster energy dissipation. This behavior is consistent with the physical interpretation of the fractional derivatives, where the parameters \(\alpha\) and \(\beta\) control the memory and damping characteristics of the material.

\section{Conclustion}
We have developed a high-order, energy-decay numerical scheme for the Maxwell equations in H-N dispersive media by overcoming a key limitation of conventional CQ methods. The analysis established that directly employing the generating function from linear multistep methods cannot produce completely monotonic sequences of order higher than one. To circumvent this barrier, we introduced a novel reconstruction of the generating function, ingeniously incorporating the properties of Pick functions. This new framework enabled the construction of a second-order, completely monotonic sequence, which is crucial for preserving the energy decay structure at the discrete level. The resulting numerical scheme is proven to be unconditionally stable. Numerical tests confirm the theoretical second-order accuracy in time and the desired monotonic energy decay, demonstrating the effectiveness of our Pick function-based approach for simulating complex dispersive electromagnetic systems.

\section*{Acknowledgements}
This work is supported by the National Natural Science Foundation of China (No. 12201322 to B.Y., No. 12401530 to G.Y., No. 12461080 to Y.L and No. 12561068 to H.L.), 
Natural Science Foundation of Inner Mongolia (No. 2025MS01003 to B.Y.),
Program for Innovative Research Team in Universities of Inner Mongolia Autonomous Region (No. NMGIRT2413 to Y.L.),
Key Project of Natural Science Foundation of Inner Mongolia Autonomous Region (No. 2025ZD036 to H.L.).
 \section*{Declarations}
\textbf{Conflict of interest} ~ The authors declare that they have no conflict of interest.
\bibliography{mybibfile}

\end{document}